\documentclass[11pt,a4paper]{article}
\usepackage{amsfonts}
\usepackage{amsmath}
\usepackage{latexsym}
\usepackage{mathrsfs}
 \newtheorem{theorem}{Theorem}[section]
 \newtheorem{remark}{Remark}[section]
 
 \newtheorem{lemma}{Lemma}[section]

\def\qed{\hbox to 0pt{}\hfill$\rlap{$\sqcap$}\sqcup$}
\numberwithin{equation}{section} 

\def\whitebox{{\hbox{\hskip 1 pt
     \vrule height 6pt depth 1.5pt
     \lower 1.5pt\vbox to 7.5pt{\hrule width
     3.2pt\vfill\hrule width 3.2pt}%
     \vrule height 6pt depth 1.5pt
     \hskip1pt}}}
\def\qed{\ifhmode\allowbreak\else\nobreak\fi\hfill\quad\nobreak
\whitebox\medbreak}
\begin{document}
\title{\Large \bf Positive Solutions for Impulsive Integral Boundary Value Problems on an Infinite Interval}

\author{Ilkay Yaslan Karaca and Sezgi Aksoy}

\date{}
\maketitle
\begin{center}
	Department of Mathematics, Ege University, 35100 Bornova, Izmir, Turkey\\
	E-mail: ilkay.karaca@ege.edu.tr, sezgiaksoy6@gmail.com

\end{center}

{\noindent \bf Abstract}\\
 \indent  This paper is concerned with the existence of positive solutions of second-order impulsive differential equations with integral boundary conditions on an infinite interval. As an application, an example is given to demonstrate our main results.

\vspace{.5cm} {\noindent \bf 2010 Mathematics subject
Classification :} 34B37, 34B15, 34B18, 34B40.

{\noindent \bf Keywords:} Impulsive differential equations, positive solutions, fixed point theorems in cones, impulsive boundary value problems, integral boundary conditions.\\
\section{\large\bf Introduction}
\label{sec1}\indent\indent 
Consider the following second-order impulsive integral boundary value problem (IBVP) with integral boundary conditions,
\begin{eqnarray} \label{a1}
\begin{cases}
\displaystyle \frac{1}{p(t)}(p(t)x'(t))'+f(t,x(t),x'(t))=0,\ \ \ \ \forall t\in J^{'}_{+} \\
\Delta x |_{t= t_k} = I_{k} (x(t_{k})),\ \ \ \ \ \ \ k=1,2,...,n  \\  
\Delta x'|_{t= t_k} =-\overline I_{k}(x(t_{k})),\ \ \ \ k=1,2,...,n  \\  
a_{1}x(0)-b_{1}\displaystyle\lim_{t\rightarrow 0^+}p(t)x'(t)=      \displaystyle\int_{0}^{\infty}g_{1}(x(s)) \psi (s)ds , \\
a_{2}\displaystyle \lim_{t\rightarrow \infty}x(t)+b_{2} \displaystyle\lim_{t\rightarrow \infty}p(t)x'(t)=\displaystyle\int_{0}^{\infty}g_{2}(x(s)\psi(s)ds, \\
\end{cases}
\end{eqnarray}
where $ J=[0,\infty)$, $J_{+}=(0, \infty)$, $J'_{+}=J_{+}\backslash\{t_{1},...,t_{n}\}$, $0<t_{1}<t_{2}<...<t_{n}$, and note $J_{0}=[0,t_{1})$, $J_{i}=(t_{i},t_{i+1}], (i=1,2,...,n)$, $\Delta x |_{t= t_k}$ and $\Delta x'|_{t= t_k}$ denote the jump of $x(t)$ and $ x^{'}(t)$ at $t=t_{k}$, i.e.,
$$\Delta x |_{t= t_k}=x(t^+_k)-x(t^-_k),\ \ \ \ \ \Delta x'|_{t= t_k}=x'(t^+_k)-x'(t^-_k),$$
where $x(t^{+}_{k})$, $x^{'}(t^{+}_{k})$ and  $x(t^{-}_{k})$, $x^{'}(t^{-}_{k})$ respent the right-hand limit and left-hand limit of $x(t)$ and $x^{'}(t)$ at $t=t_{k},\ k=1,2,...,n$, respectively.

Throughout this paper, we assume that the folllowing conditions hold;
\begin{itemize}
	\item[$(H1)$] \label{h1} $a_{1}, a_{2}$, $b_{1},b_{2}\in J$ with $D=a_{2}b_{1}+a_{1}b_{2}+a_{1}a_{2}B(0,\infty)>0 $ in which $ B(t,s)=\displaystyle\int_{t}^{s}\frac{d\sigma}{p(\sigma)}$.
	\item[$(H2)$] \label{h2} $f\in C(J_{+}\times J_{+}\times\mathbb{R}, J_{+})$ and also,
	\begin{eqnarray*}
		f(t,x,y)\leq k(t)h(x,y),\ \ \ \ \ (t\in J_{+})
	\end{eqnarray*}
	where $h\in C(J_{+}\times\mathbb{R},J_{+})$ and for $t\in J_{+}$, $x,\ y$ in a bounded set, $h(x,y)$ is bounded and $k\in C(J_{+},J_{+})$ .
	\item[$(H3)$] $g_{1},g_{2} : J_{+}\rightarrow J_{+}$ are continuous, nondecreasing functions, and for $t\in J_{+}$, $x$ in a bounded set, $g_{1}(x),g_{2}(x)$ are bounded.
	\item[$(H4)$] $I_{k}, \overline I_{k} \in C(J_{+},J_{+})$ are bounded functions where
	\begin{eqnarray*}
		[b_{2}+a_{2}B(t_{k},\infty)]\overline I_{k}(x(t_{k}))- \displaystyle\frac{a_{2}}{p(t_{k})}I_{k}(x(t_{k}))>0,\ \ \ (k=1,2,...,n).
	\end{eqnarray*} 
	\label{h5}\item[$(H5)$] $\psi : J\rightarrow J$ is a continuous function with $\displaystyle\int_{0}^{\infty}\psi(s)ds<\infty$.
	\label{h6}\item[$(H6)$] $p\in C(J,J_{+})\cap C^1(J_{+},J_{+})$ with $p>0$ on $J_{+}$, and $\displaystyle\int_{0}^{\infty}\frac{ds}{p(s)}<\infty$.
\end{itemize}

Impulsive differential equations encountered in physics, chemical technology, population dynamics, biotecnology, economics etc. (see \cite{book.1} and the references there in) have become more important in recent years due to the appearance of some important in recent years due to the appearence of some mathematical models of the actual processes. A significant development was observed in theory of impulsive differential equations with fixed time of pulses; see the monographs by Bainov and Simeorov \cite{book.2}, Lakshmikantham, et al. \cite{book.3}, Samoilenko and Prestyuk \cite{book.4}, Benchohra, et al. \cite{book.5} and the papers \cite{article.24,article.26,article.13,article.14,article.10,article.9,article.25,article.23}.

The existence and multiplicity of positive solutions for linear and nonlinear second-order impulsive dynamic equations have been extensively studied, see \cite{article.4,article.11,article.5,article.20,article.6}. Due to the fact that an infinite interval is noncompact, the discussion about boundary value problems on the half-line more complicated, in particular, for impulsive IBVP on an infinite interval, few works were done, see \cite{article.7,article.8}. There is not work on positive solutions for double impulsive IBVP on an infinite interval expect that in \cite{article.22,article.3}.

In \cite{article.3}, Zhang, Yang and Feng studied the following double impulsive IBVP:
\begin{eqnarray*}
	\begin{cases}
		-x^{''}(t)= f(t,x(t),x^{'}(t)),\ \ \ t\in J,\ \ \ t\neq t_{k},\\
		\Delta x |_{t= t_k}=I_{k} (x(t_{k})),\ \ \ \ \ \ \ k=1,2,...\\
		\Delta x^{'} |_{t= t_k}=\overline I_{k}(x(t_{k})),\ \ \ \ \ \ k=1,2,...\\ x(0)=\displaystyle\int_{0}^{\infty} g(t)x(t)dt,\ \ \ \ \ x^{'}(\infty)=0.
	\end{cases} 
\end{eqnarray*}
Using the fixed point theorem in cones, they obtained criteria for existence of the multiple positive solutions.

In \cite{article.22}, Yu, Wang and Guo discussed the existence and multiple positive solutions for the following nonlinear second-order double impulsive integral boundary value problems:
\begin{eqnarray*}
	\begin{cases}
		(\phi_{p}(x^{'}(t)))^{'}+a(t)f(t,x(t),x^{'}(t))=0,\ \ t\in J,\ \ t\neq t_{k},\\
		\Delta x |_{t= t_k}=I_{k} (x(t_{k})),\ \ \ \ \ \ \ \ \ \ k=1,2,...\\
		\Delta \phi_{p}(x^{'})|_{t= t_k}=\overline I_{k}(x(t_{k})),\ \ \ \ k=1,2,...\\ x(0)=\displaystyle\int_{0}^{\infty} g(t)x(t)dt,\ \ \ \ \ x^{'}(\infty)=0.
	\end{cases} 
\end{eqnarray*}

Motivated by the above works, in this study, we consider the existence of two positive solutions for the second-order double impulsive integral boundary value problem \eqref{a1}. Our boundary conditions are more general. Hence, these results can be considered as a contribution to this field.

The present paper is organized as follows. In section 2, we present some preliminaries and lemmas which are key tools for our main results. We give and prove our main results in section 3. Finally, in section 4, we give an example to demonstrate our results.

\section{Preliminaries and Lemmas} \label{sec:1}

In this section, we will employ several lemmas to prove the main results in this paper.\\
Set\\

$PC(J)=\big\{x:J\rightarrow\mathbb{R}:\ x\in C(J^{'})$,  $\ x(t^{+}_{k})$ and $x(t^{-}_{k})$ exist and $ x(t^{-}_{k})=x(t_{k}),$ $\ \ \ \ 1\leq k\leq n \big\}$.\\
$PC^{1}(J)=\big\{x\in PC(J):\ x^{'} \in C(J^{'}),\ x^{'}(t^{+}_{k})$ and $x^{'}(t^{-}_{k})$ exist and $x^{'}(t^{-}_{k})=x^{'}(t_{k})\big\}$.\\

$BPC^{1}(J)=\big\{x\in PC^{1}(J):\ \displaystyle\lim_{t\rightarrow \infty} x(t)$ exists, and $\ \displaystyle\sup_{t\in J} |x^{'}(t)|<\infty \big\}$.\\ \\
It is easy to see that $BPC^{1}(J)$ is a Banach space with the norm
$$\big\|x\big\|= \displaystyle \sup_{t\in J} \big\{|x(t)|+|x^{'}(t)|\big\}.$$
A function $ x\in PC^{1}(J)\cap C^{2}(J^{'}_{+})$ is called a positive solution of the impulsive IBVP $\eqref{a1}$ if $x(t)>0$ for all $t\in J$ and $x(t)$ satisfies $\eqref{a1}$.\\ \\
We define a cone $K\subset BPC^{1}(J)$ as follows:
$$K=\big\{x\in BPC^{1}(J) : x(t)>0,\ t \in J_{+}\big\}.$$
K is a positive cone in $BPC^{1}(J)$.

By $\theta$ and $\varphi$  we denote the solutions of the corresponding homogeneous equation
\begin{eqnarray} \label{b2}
\frac{1}{p(t)}(p(t)x'(t))'=0 , \ t\in (0,\infty),
\end{eqnarray}
under the initial conditions,
\begin{eqnarray}\label{b3}
\nonumber \theta(0)=b_{1},\ \ \displaystyle\lim_{t\rightarrow 0^+} p(t)\theta'(t)=a_{1}, \\ \\ \nonumber
\displaystyle\lim_{t\rightarrow\infty}\varphi(t)=b_{2},\ \displaystyle\lim_{t\rightarrow\infty} p(t)\varphi'(t)=-a_{2}. 
\end{eqnarray}

Using the initial conditions $\eqref{b3}$, we can deduce from equation $\eqref{b2}$ for $\theta(t)$ and $\varphi(t)$, the following equations:
\begin{eqnarray}
\label{b4}\theta(t)=b_{1}+a_{1}\displaystyle\int_{0}^{t}\frac{ds}{p(s)}, \\
\label{b5}\varphi(t)=b_{2}+a_{2}\displaystyle\int_{t}^{\infty}\frac{ds}{p(s)}.
\end{eqnarray}

Let $G(t,s)$ be the Green Function for $\eqref{a1}$ is given by
\begin{eqnarray}\label{b6}
G(t,s)=\frac{1}{D}
\begin{cases}
\theta(t)\varphi(s), \ \ 0 \leq t \leq s<\infty, \\
\theta(s)\varphi(t), \ \ 0 \leq s \leq t<\infty. \\ 
\end{cases}
\end{eqnarray}
where $\theta(t)$ and $\varphi(t)$ are given in $\eqref{b4}$ and $\eqref{b5}$ respectively. 

\begin{lemma}
	Suppose that $(H1)-(H6)$ are satisfied. Then $x\in PC^{1}(J)\cap C^{2}(J^{'}_{+})$ is a solution of the impulsive IBVP $\eqref{a1}$ if and only if $x(t)$ is a solution of the following integral equation
	\begin{eqnarray*}
		\nonumber&x(t)=&\displaystyle\int_{0}^{\infty}G(t,s)p(s)f(s,x(s),x^{'}(s))ds
		+\frac{\varphi(t)}{D}\displaystyle\int_{0}^{\infty}g_{1}(x(s))\psi(s)ds\\&
		&
		+\frac{\theta(t)}{D} \displaystyle\int_{0}^{\infty}g_{2}(x(s))\psi(s)ds
		+\displaystyle\sum_{k=1}^{n}G(t,t_{k})\overline I_{k}(x(t_{k}))
		+\displaystyle\sum_{k=1}^{n}p(t_{k})G_{s}(t,s)|_{s=t_{k}}I_{k}(x(t_{k})),
	\end{eqnarray*}
	where G(t,s) is given by $\eqref{b6}$.
\end{lemma}
\begin{remark}
	It is easy to prove the following properties of G(t,s),
	\begin{itemize}
		\item[$(1)$] $ G(t,s) $ is continuous on $J_{+}\times J_{+}$,
		\item[$(2)$] For each $s\in J_{+}$, G(t,s) is continuously differentiable on $J_{+}$ except $t=s$,
		\item[$(3)$] $\displaystyle\frac{\partial G(t,s)}{\partial t} \big|_{t=s^+}- \displaystyle\frac{\partial G(t,s)}{\partial t} \big|_{t=s^-}=\frac{1}{p(s)}$,
		\item[$(4)$] $G(t,s)\leq G(s,s)<\infty$, and $G_{s}(t,s)\leq G_{s}(t,s)\big|_{t=s}<\infty$,
		\item[$(5)$] $|G_{t}(t,s)|\leq \displaystyle\frac{c}{p(t)}G(s,s),$ and $|G_{st}(t,s)|\leq \displaystyle\frac{c}{p(t)}G_{s}(t,s)\big|_{t=s},$\\
		where
		\begin{eqnarray}\label{b7}
		c=\frac{\max\{a_{1},a_{2}\}}{\min\{b_{1},b_{2}\}},
		\end{eqnarray} 
		\item[$(6)$] $\overline{G}(s)=\displaystyle\lim_{t\rightarrow \infty} G(t,s)=\displaystyle\frac{b_{2}}{D} \ \theta(s) \leq G(s,s)<\infty,$
		\item[$(7)$] $\overline G^{'}(s)=\displaystyle\lim_{t\rightarrow\infty}G_{s}(t,s)= \frac{b_{2}}{D}\ \theta^{'}(s)\leq G_{s}(t,s)\big|_{t=s}<\infty,$
		\item[$(8)$] For any $t\in[a,b]\subset(0,\infty)$ and $s\in [0,\infty),$ we have \\
		$$G(t,s) \geq wG(s,s),\ \ and \ \ \ G_{s}(t,s)\geq wG_{s}(t,s)\big|_{t=s},$$\\
		where
		\begin{eqnarray}\label{b8}
		w= \min \big\{\frac{b_{1}+a_{1}B(0,a)}{b_{1}+a_{1}B(0,\infty)}, \frac{b_{2}+a_{2}B(b,\infty)}{b_{2}+a_{2}B(0,\infty)} \big\}.
		\end{eqnarray}
		Obviously, $ 0<w<1$.
	\end{itemize}
\end{remark}
Define
\begin{eqnarray}\label{b9}
& \nonumber (Tx)(t)=&\displaystyle\int_{0}^{\infty}G(t,s)p(s)f(s,x(s),x'(s))ds
+\frac{\varphi(t)}{D}\displaystyle\int_{0}^{\infty}g_{1}(x(s))\psi(s)ds \\&
&\ \ 
+\frac{\theta(t)}{D} \displaystyle\int_{0}^{\infty}g_{2}(x(s))\psi(s)ds+\displaystyle\sum_{k=1}^{n}G(t,t_{k})\overline I_{k}(x(t_{k})) \\& \nonumber
&\ \
+\displaystyle\sum_{k=1}^{n}p(t_{k})G_{s}(t,s)|_{s=t_{k}}I_{k}(x(t_{k})),
\end{eqnarray}
where $G$ is defined by as in $\eqref{b6}$.

Obviously, the impulsive IBVP $\eqref{a1}$ has a solution $x$ if and only if $x\in K$ is a fixed point of the operator $T$ defined by $\eqref{b9}$.

It is convenient to list the following condition which is to be used in our theorems:
\begin{itemize}
	\label{h7} \item[$(H7)$] $0<\displaystyle\int_{0}^{\infty}G(s,s)p(s)k(s)ds<\infty$.
\end{itemize}
As we know that the Ascoli-Arzela Theorem does not hold in infinite interval $J$, we need the following compactness critarion:

\begin{lemma}(\cite{article.15})
	Let $M\subset BPC^{1}(J)$. Then M is relatively compact in $BPC^{1}(J)$ if the following conditions hold.
	\begin{itemize}
		\item[$(i)$] M is uniformly bounded in $BPC^{1}(J)$.
		\item[$(ii)$] The function belonging to M are equicontinuous on any compact interval of $[0,\infty)$.
		\item[$(iii)$] The functions from M are equiconvergent, that is, for any given $\varepsilon>0$, there exist a \\ $T=T(\varepsilon)>0$ such that $ |f(t)-f(\infty)|<\varepsilon$ for any $t>T,\ f\in M$.
	\end{itemize}
	The main tool of this work is a fixed point theorem in cones.
\end{lemma}

\begin{lemma}(\cite{article.16})
	Let X be an Banach space and K is a positive cone in X. Assume that $\Omega_{1},\Omega_{2}$ are open subsets of X with $0\in \Omega_{1},\ \overline\Omega_{1}\subset \Omega_{2}$. Let $T:K\cap (\overline\Omega_{2}\backslash\Omega_{1})\rightarrow K$ be a completely continuous operator such that
	\begin{itemize}
		\item[$(i)$] $\big\|Tx\big\|\leq\big\|x\big\|$ for all $x\in K\cap\partial\Omega_{1}$.
		\item[$(ii)$] There exists a $\Phi\in K$ such that $x\neq Tx+\lambda\Phi$, for all $x\in K\cap\partial\Omega_{2}$ and $\lambda>0$.
	\end{itemize}
	Then T has a fixed point in $K\cap (\overline\Omega_{2}\backslash\Omega_{1})$.
\end{lemma}

\begin{lemma}
	If (H1)-(H7) are satisfied, then for any bounded open set $\Omega\subset BPC^{1}(J)$, $T:K\cap\overline\Omega\rightarrow K$ is a completely continuous operator.
\end{lemma}
\noindent{\em Proof}. For any bounded open set $\Omega\subset BPC^{1}(J)$, there exists a constant $M>0$ such that $\big\|x\big\|\leq M$ for any $x\in\overline\Omega$.

First, we show $T:K\cap\overline\Omega\rightarrow K$ is well defined. Let $x\in K\cap\overline\Omega $. From $(H2)$, $(H3)$ and $(H4)$, we have
\begin{eqnarray}\label{b10}
S_{M}&=	sup \{S_{1},S_{2},S_{3},S_{4},S_{5} \},
\end{eqnarray}
where
\begin{eqnarray*}
	S_{1}&=	&sup \{ h(x,y): |x|+|y| \leq M\}<\infty, \\
	S_{2}&=	&sup \{ I_{k}(x) :\ 0\leq x\leq M \}, \\
	S_{3}&=	&sup \{ \overline I_{k}(x) :\ 0\leq x\leq M \}, \\
	S_{4}&=	&sup \{ g_{1}(x) :\ 0\leq x\leq M \}, \\
	S_{5}&=	&sup \{ g_{2}(x) :\ 0\leq x\leq M \}.
\end{eqnarray*} \\

Let $t_{1},t_{2}\in J,\ t_{1}<t_{2},\ $ then
\begin{eqnarray}\label{b11}
\displaystyle\int_{0}^{\infty}\big|G(t_{1},s)-G(t_{2},s)\big|p(s)k(s)ds\leq 2\displaystyle\int_{0}^{\infty}G(s,s)p(s)k(s)ds <\infty .
\end{eqnarray}
Hence, by the Lebesgue dominated convergence theorem, we have for any $t_{1},t_{2}\in J$, $x\in K\cap\overline\Omega$, and the fact that $G(t,s)$ is continuous, we have

\begin{eqnarray}\label{b12}
&\big|(Tx)(t_{1})-(Tx)(t_{2})\big| \nonumber &\leq \displaystyle\int_{0}^{\infty}\big|G(t_{1},s)-G(t_{2},s)\big|p(s)f(s,x(s),x^{'}(s))ds \\&
& \ \ \nonumber +\frac{|\varphi(t_{1})-\varphi(t_{2})|}{D} \displaystyle\int_{0}^{\infty}g_{1}(x(s))\psi(s)ds \\&
& \ \ \nonumber +\frac{|\theta(t_{1})-\theta(t_{2})|}{D}  \displaystyle\int_{0}^{\infty}g_{2}(x(s))\psi(s)ds \\&
&\ \ \nonumber +\displaystyle\sum_{k=1}^{n}\big|G(t_{1},t_{k})-G(t_{2},t_{k})\big|\overline I_{k}(x(t_{k})) \\&
&\ \ \nonumber
+\displaystyle\sum_{k=1}^{n}\big|p(t_{k})G_{s}(t,s)|_{\substack{t=t_{1}\\s=t_{k}}}-p(t_{k})G_{s}(t,s)|_{\substack{t=t_{2}\\s=t_{k}}}\big|I_{k}(x(t_{k})) \\&
&\nonumber
\leq S_{M} \bigg\{\displaystyle\int_{0}^{\infty}\big|G(t_{1},s)-G(t_{2},s)\big|p(s)k(s)ds \\&
&\ \ \nonumber+\bigg[\frac{|\varphi(t_{1})-\varphi(t_{2})|}{D}+\frac{|\theta(t_{1})-\theta(t_{2})|}{D}\bigg]\displaystyle\int_{0}^{\infty}\psi(s)ds 
\\&
&\ \ \nonumber +\displaystyle\sum_{k=1}^{n}\big|G(t_{1},t_{k})-G(t_{2},t_{k})\big|\\&
&\ \ \nonumber +\displaystyle\frac{1}{D}\displaystyle\sum_{t_{k}\leq t_{1}}p(t_{k})\theta^{'}(t_{k})\big|\varphi(t_{1})-\varphi(t_{2})\big| \bigg\} \\&
&\ \ \nonumber +\displaystyle\frac{1}{D}\displaystyle\sum_{t_{2}\leq t_{k}}p(t_{k})\big|\varphi^{'}(t_{k})\big|\big|\theta(t_{1})-\theta(t_{2})\big| \\&
&\ \ \nonumber +\displaystyle\frac{1}{D}\displaystyle\sum_{t_{1}\leq t_{k}\leq t_{2}}p(t_{k})\big|\theta^{'}(t_{1})\varphi^{'}(t_{k})-\theta^{'}(t_{k})\varphi^{'}(t_{2})\big| \bigg\} \\&
& \rightarrow 0 \ \ \ as \ \ t_{1}\rightarrow t_{2},
\end{eqnarray}

\begin{eqnarray}\label{b13}
&\big|(Tx)^{'}(t_{1})-(Tx)^{'}(t_{2})\big| &\leq \nonumber
S_{M} \bigg\{ \frac{a_{2}}{D}\bigg|\frac{1}{p(t_{1})}-\frac{1}{p(t_{2})}\bigg|\displaystyle\int_{0}^{t_{1}}
\theta(s)p(s)k(s)ds \\&
&\ \ \nonumber\ \ \ +\frac{a_{1}}{Dp(t_{1})}\displaystyle\int_{t_{1}}^{t_{2}}\varphi(s)p(s)k(s)ds \\&
&\ \ \nonumber\ \ \ +\frac{a_{1}}{D}\bigg|\frac{1}{p(t_{1})}-\frac{1}{p(t_{2})}\bigg|\displaystyle\int_{t_{2}}^{\infty}\varphi(s)p(s)k(s)ds \\&
&\ \ \nonumber\ \ \
+\frac{a_{2}}{Dp(t_{2})} \displaystyle\int_{t_{1}}^{t_{2}}\theta(s)p(s)k(s)ds\\&
&\ \ \nonumber\ \ \ +\frac{a_{2}}{D}\bigg|\frac{1}{p(t_{2})}-\frac{1}{p(t_{1})}\bigg|\displaystyle\int_{0}^{\infty}\psi(s)ds
 \\&
&\ \ \nonumber\ \ \ +\frac{a_{1}}{D}\bigg|\frac{1}{p(t_{1})}-\frac{1}{p(t_{2})}\bigg|\displaystyle\int_{0}^{\infty}\psi(s)ds
\end{eqnarray}

\begin{eqnarray}  
& &\ \ \nonumber\ \ \ +\displaystyle\frac{a_{1}}{D}\bigg|\frac{1}{p(t_{1})}-\frac{1}{p(t_{2})}\bigg|\displaystyle\sum_{t_{2}\leq t_{k}}\big|\varphi(t_{k})+p(t_{k})\varphi^{'}(t_{k})\big|
\\&
&\ \ \nonumber\ \ \ +\frac{a_{2}}{D}\bigg|\frac{1}{p(t_{1})}-\frac{1}{p(t_{2})}\bigg| \displaystyle\sum_{t_{k}\leq t_{1}}\big|\theta(t_{k})+p(t_{k})\theta^{'}(t_{k})\big| \\&
&\ \ \nonumber\ \ \ 
+\frac{a_{1}}{Dp(t_{1})}\displaystyle\sum_{t_{1}\leq t_{k}\leq t_{2}}\big|\varphi(t_{k})+p(t_{k})\varphi^{'}(t_{k})\big| \\&
&\ \ \nonumber\ \ \ +\frac{a_{2}}{Dp(t_{2})}\displaystyle\sum_{t_{1}\leq t_{k}\leq t_{2}} \big|\theta(t_{k})+p(t_{k})\theta^{'}(t_{k})\big| \bigg\} \\&
&\ \ \rightarrow 0 \ \ \ as \ \ t_{1}\rightarrow t_{2}.
\end{eqnarray}
Thus, $Tx \in PC^{1}(J) $.
We can show that $Tx\in BPC^{1}(J)$.
Then by $(H5)$, $(H7)$, the properties $(5)$, $(6)$, $(7)$ of Remark $1$ and the Lebegue dominated convergence theorem, we have
\begin{eqnarray}\label{b14}
&\displaystyle\lim_{t\rightarrow \infty}(Tx)(t)&= \nonumber\displaystyle\int_{0}^{\infty} \overline G(s)p(s)f(s,x(s),x^{'}(s))ds
+\frac{\varphi(\infty)}{D}\displaystyle\int_{0}^{\infty}g_{1}(x(s))\psi(s)ds\\\nonumber &
&\ \ \ \ +\frac{\theta(\infty)}{D} \displaystyle\int_{0}^{\infty}g_{2}(x(s))\psi(s)ds +\displaystyle\sum_{k=1}^{n}\overline G(t_{k})\overline I_{k}(x(t_{k}))
\\\nonumber &
&\ \ \ \ 
+\displaystyle\sum_{k=1}^{n}p(t_{k})\overline G^{'}(t_{k})I_{k}(x(t_{k})) \\&
&<\infty
\end{eqnarray}
and
\begin{eqnarray}\label{b15}
&\big|(Tx)^{'}(t)\big|\nonumber &\leq S_{M} \bigg\{ \frac{c}{p(t)}\displaystyle\int_{0}^{\infty}G(s,s)p(s)k(s)ds+ \frac
{\max\{a_{1},a_{2}\}}{Dp(t)}\displaystyle\int_{0}^{\infty}\psi(s)ds \hspace{3cm}\\&
& \ \ \ \ \nonumber +\frac{c}{p(t)} \displaystyle\sum_{k=1}^{n}G(t_{k},t_{k}) +\frac{c}{p(t)}\displaystyle\sum_{k=1}^{n}p(t_{k})G_{s}(t,s)|_{\substack{t=t_{k}\\s=t_{k}}} \bigg\} \\ &
&<\infty.
\end{eqnarray}
Therefore, $\displaystyle\sup_{t\in J}|(Tx)^{'}(t)|<\infty$.\\

Hence $T:K\cap\overline\Omega\rightarrow K $ is well defined. \\

Next, we prove that $T$ is continuous. Let $ x_{n}\rightarrow x$ in $K\cap\overline\Omega$, then $ \|x_{n}\| \leq M \ \ (n=1,2,...)$. We will show that $Tx_{n}\rightarrow Tx$. For any $\varepsilon>0$, by $(H7)$ there exists a constant $A_{0}>0$ such that
\begin{eqnarray}\label{b16}
S_{M}\displaystyle\int_{A_{0}}^{\infty}G(s,s)p(s)k(s)ds\leq \frac{\varepsilon}{12}.
\end{eqnarray} 

On the other hand, by the continuity of $f(t,u,v)$ on $(0,A_{0}]\times J_{+} \times \mathbb{R}$, the continuities of $g_{1}$, $g_{2}$ on $J_{+}$ and the continuities of $I_{k}$, $\overline I_{k}$ on $J_{+}$, for the above $\varepsilon>0$, there exists a $\delta>0$ such that, for any $u$, $v$, $u_{1}$, $v_{1}$, satisfying $|u|+|v|<M$, and $|u_{1}|+|v_{1}|<M$, $\big|u-u_{1}\big|+\big|v-v_{1}\big|<\delta$,
\begin{eqnarray}\label{b17b}
\nonumber& &\big|f(s,u,v)-f(s,u_{1},v_{1})\big|< \frac{\varepsilon}{6}\bigg(\displaystyle\int_{0}^{A_{0}}G(s,s)p(s)ds\bigg)^{-1},\\ \nonumber
\nonumber& &\big|g_{1}(u)-g_{1}(u_{1})\big|<\frac{\varepsilon}{6}\bigg( \frac{\varphi(0)}{D}\displaystyle\int_{0}^{\infty}\psi(s)ds \bigg)^{-1},\\
& &	\big|g_{2}(u)-g_{2}(u_{1})\big|<\frac{\varepsilon}{6}\bigg( \frac{\theta(\infty)}{D}\displaystyle\int_{0}^{\infty}\psi(s)ds \bigg)^{-1},\\\nonumber 
& &\big|{I_{k}}(u(t_{k}))-{I_{k}}(u_{1}(t_{k})) \big| <\frac{\varepsilon}{6}\bigg( \displaystyle\sum_{k=1}^{n}p(t_{k})G_{s}(t,s)|_{\substack{t=t_{k}\\s=t_{k}}} \bigg)^{-1},\\\nonumber
& &\big|\overline I_{k}(u(t_{k}))-\overline I_{k}(u_{1}(t_{k})) \big| <\frac{\varepsilon}{6}\bigg( \displaystyle\sum_{k=1}^{n}G(t_{k},t_{k}) \bigg)^{-1}.
\end{eqnarray}

From the fact that $\|x_{n}-x\|\rightarrow 0$ as $n\rightarrow\infty$, for above $\delta$, there exists a sufficiently large number N such that, when $n>N$, we have, for $t\in(0,A_{0}]$,
\begin{eqnarray}
\hspace{1cm} |x_{n}(t)-x(t)|+|x^{'}_{n}(t)-x^{'}(t)| \leq \|x_{n}-x\|<\delta.
\end{eqnarray}

By $\eqref{b16}$-$\eqref{b17b}$, we have, for $n>N$,
\begin{eqnarray*}
	&\big|(Tx_{n})(t)-(Tx)(t)\big| &\leq  \bigg| \displaystyle\int_{0}^{\infty}G(s,s)p(s)\big[f(s,x_{n}(s),x^{'}_{n}(s))-f(s,x(s),x^{'}(s))\big]ds \\&	& \ \ \
	+\frac{\varphi(0)}{D}\displaystyle\int_{0}^{\infty}\big[g_{1}(x_{n}(s))-g_{1}(x(s))\big]\psi(s)ds \\&	& \ \ \ +\frac{\theta(\infty)}{D}\displaystyle\int_{0}^{\infty}\big[g_{2}(x_{n}(s))-g_{2}(x(s))\big]\psi(s)ds \\&	&\ \ \ 
	+\displaystyle\sum_{k=1}^{n}G(t_{k},t_{k}) \big[ \overline I_{k}(x_{n}(t_{k}))-\overline I_{k}(x(t_{k}))\big] \\&	&\ \ \
	+\displaystyle\sum_{k=1}^{n}p(t_{k})G_{s}(t,s)|_{\substack{t=t_{k}\\s=t_{k}}}\big[I_{k}(x_{n}(t_{k}))-I_{k}(x(t_{k}))\big]\bigg|  \\&	&
	\leq  \displaystyle\int_{0}^{A_{0}}G(s,s)p(s) \big|f(s,x_{n}(s),y_{n}(s))-f(s,x(s),y(s))\big|ds \\& & \ \ \ +2S_{M}\displaystyle\int_{A_{0}}^{\infty}G(s,s)p(s)k(s)ds \\&	& \ \ \ +\frac{\varphi(0)}{D}\displaystyle\int_{0}^{\infty}\big|g_{1}(x_{n}(s))-g_{1}(x(s)) \big|\psi(s)ds \\&	& \ \ \ +\frac{\theta(\infty)}{D}\displaystyle\int_{0}^{\infty}\big|g_{2}(x_{n}(s))-g_{2}(x(s))\big|\psi(s)ds \\&	& \ \ \
	+\displaystyle\sum_{k=1}^{n}G(t_{k},t_{k}) \big|\overline I_{k}(x_{n}(t_{k}))-\overline I_{k}(x(t_{k}))\big| \\&	& \ \ \
	+\displaystyle\sum_{k=1}^{n}p(t_{k})G_{s}(t,s)|_{\substack{t=t_{k}\\s=t_{k}}}\big|I_{k}(x_{n}(t_{k}))-I_{k}(x(t_{k}))\big| \\&	&
	< \frac{\varepsilon}{6} +\frac{\varepsilon}{6}+\frac{\varepsilon}{6}+\frac{\varepsilon}{6}+\frac{\varepsilon}{6}+\frac{\varepsilon}{6} =\varepsilon.
\end{eqnarray*}
Similarly, we can see that when $\|x_{n}-x\| \rightarrow 0$ as $n\rightarrow \infty$, $|(Tx_{n})^{'}(t)-(Tx)^{'}(t)| \rightarrow 0$ as $n\rightarrow \infty$. This implies that $T$ is a continuous operator.\\

Finally we show that $T:K\cap\overline\Omega \rightarrow K $ is a compact operator. In fact for any bounded set $D\subset\overline\Omega$, there exists a constant $R>0$ such that $\|x\|\leq R$ for any $x\in K\cap D$. Hence, we have
\begin{eqnarray*}
	|(Tx)(t)|& &\leq \bigg|\displaystyle\int_{0}^{\infty}G(s,s)p(s)f(s,x(s),x^{'}(s))ds
	+\frac{\varphi(0)}{D}\displaystyle\int_{0}^{\infty}g_{1}(x(s))\psi(s)ds
	\\& &\ \ \
	+\frac{\theta(\infty)}{D} \displaystyle\int_{0}^{\infty}g_{2}(x(s))\psi(s)ds
\end{eqnarray*}

\begin{eqnarray*}
	&	&\ \ \ +\displaystyle\sum_{k=1}^{n}G(t_{k},t_{k})\overline I_{k}(x(t_{k}))
	+\displaystyle\sum_{k=1}^{n}p(t_{k})G_{s}(t,s)|_{\substack{t=t_{k}\\s=t_{k}}} I_{k}(x(t_{k}))\bigg|\\&
	& \leq S_{R} \bigg( \displaystyle\int_{0}^{\infty}G(s,s)p(s)k(s)ds
	+\frac{\varphi(0)}{D}\displaystyle\int_{0}^{\infty}\psi(s)ds\\&
	&\ \ \ 
	+\frac{\theta(\infty)}{D} \displaystyle\int_{0}^{\infty}\psi(s)ds +\displaystyle\sum_{k=1}^{n}G(t_{k},t_{k})
	+\displaystyle\sum_{k=1}^{n}p(t_{k})G_{s}(t,s)|_{\substack{t=t_{k}\\s=t_{k}}} \bigg) \\&
	&\ < \infty.
\end{eqnarray*}
From $\eqref{b15}$, we get $|(Tx)^{'}(t)|<\infty$, for $t\in J$. \\
Therefore, $T(K\cap D)$ is uniformly bounded in $BPC^{1}(J)$. \\

Given $r>0$, for any $x\in K\cap D$, as the proof of $\eqref{b12}, \eqref{b13}$, we can get that $F=\{Tx:x\in K\cap D\} $ are equicontinuous on $[0,r]$. Since $r>0$ arbitrary, $F$ is locally equicontinuous on $J_{+}$. By $(H5)$, $(H7)$ the properties $(5)$, $(6)$, $(7)$ and the Lebesgue dominated converges theorem, we get
\begin{eqnarray*}
	&|(Tx)(t)-(Tx)(\infty)|
	& \leq S_{R} \bigg( \displaystyle\int_{0}^{\infty}\big|G(t,s)-\overline G(s)\big|p(s)k(s)ds \\&
	&\ \ +\frac{|\varphi(t)-\varphi(\infty)|}{D} \displaystyle\int_{0}^{\infty}\psi(s)ds \\&
	&\ \  +\frac{|\theta(t)-\theta(\infty)|}{D}  \displaystyle\int_{0}^{\infty}\psi(s)ds  +\displaystyle\sum_{k=1}^{n}\big|G(t,t_{k})- \overline G(t_{k})\big|\\&
	&\ \  +\displaystyle\frac{1}{D}\displaystyle\sum_{t_{k}\leq t}p(t_{k})\theta^{'}(t_{k})\big|\varphi(t)-\varphi(\infty)\big| \\&
	&\ \ +\displaystyle\frac{1}{D}\displaystyle\sum_{t\leq t_{k}}p(t_{k})\big|\theta(t)\varphi^{'}(t_{k})-\theta(t_{k})\varphi^{'}(\infty)\big|\bigg) \\& \\
	& \rightarrow 0 \ \ \ as \ \ \ t \rightarrow \infty
\end{eqnarray*}
and 
\begin{eqnarray}\label{b17}
&|(Tx)^{'}(t)-(Tx)^{'}(\infty)| & \leq \nonumber
S_{R} \bigg\{ \frac{a_{2}}{D}\bigg|\frac{1}{p(t)}-\frac{1}{p(\infty)}\bigg|\displaystyle\int_{0}^{t}
\theta(s)p(s)k(s)ds \\\nonumber\\&
&\ \ \nonumber +\frac{a_{1}}{D}\bigg|\frac{1}{p(t)}-\frac{1}{p(\infty)}\bigg|\displaystyle\int_{t}^{\infty}\varphi(s)p(s)k(s)ds\\\nonumber\\&
&\ \ \nonumber +\frac{a_{2}}{D}\bigg|\frac{1}{p(t)}-\frac{1}{p(\infty)}\bigg|\displaystyle\int_{0}^{\infty}\psi(s)ds \\\nonumber\\&
&\ \ \nonumber +\frac{a_{1}}{D}\bigg|\frac{1}{p(t)}-\frac{1}{p(\infty)}\bigg|\displaystyle\int_{0}^{\infty}\psi(s)ds\\\nonumber\\&
&\ \ \nonumber +\frac{a_{1}}{Dp(t)}\displaystyle\sum_{t\leq t_{k}}\big|\varphi(t_{k})+p(t_{k})\varphi^{'}(t_{k})\big|
\end{eqnarray}

\begin{eqnarray}
\nonumber\\&
&\ \ \nonumber +\frac{a_{2}}{D}\bigg|\frac{1}{p(t)}-\frac{1}{p(\infty)}\bigg| \displaystyle\sum_{t_{k}\leq t}\big|\theta(t_{k})+p(t_{k})\theta^{'}(t_{k})\big| \\\nonumber\\&
&\ \ \nonumber  +\frac{a_{2}}{Dp(\infty)} \displaystyle\sum_{t\leq t_{k}}p(t_{k})\big|\theta(t_{k})+p(t)\theta^{'}(t_{k})\big| \bigg\} \\\nonumber\\&
&\rightarrow 0 \ \ \ as \ \ t\rightarrow \infty.
\end{eqnarray}

Hence $T(K\cap D)$ is equiconvergent in $BPC^{1}(J)$. By Lemma $2.2$, we have that $F$ is relatively compact in $BPC^{1}(J)$. Therefore, $T:K\cap\overline\Omega \rightarrow K$ is completely continuous. 
\begin{flushright}
	$\square$
\end{flushright}

\section{ Main Results }

For convenience and simplicity in the following discussion, we use following notations:

\begin{eqnarray*}
	f_{0}=\displaystyle\lim_{|x|+|y|\rightarrow 0} inf  \min_{t\in [a,b]}\frac{f(t,x,y)}{|x|+|y|},  \hspace{1cm}      f_{\infty}=\displaystyle\lim_{|x|+|y|\rightarrow \infty} inf  \min_{t\in [a,b]}\frac{f(t,x,y)}{|x|+|y|},\\ \\
	g_{i_{0}}=\displaystyle\lim_{x\rightarrow 0} inf \frac{g_{i}(x)}{x}\ \ (1\leq i\leq 2), \hspace{1cm}
	g_{i_{\infty}}=\displaystyle\lim_{x\rightarrow \infty} inf \frac{g_{i}(x)}{x}\ \ (1\leq i\leq 2),\\ \\ 
	h^{q}=\displaystyle\lim_{|x|+|y|\rightarrow q}sup \frac{h(x,y)}{|x|+|y|}, \hspace{2cm}\ \	g_{i}^{q}=\displaystyle\lim_{x\rightarrow q}sup \frac{g_{i}(x)}{x}\ \ (1\leq i\leq 2),
\end{eqnarray*}
\begin{eqnarray*}
	I_{0}(k)=\displaystyle\lim_{x\rightarrow 0} inf \frac{I_{k}(x)}{x},\hspace{4cm} I_{\infty}(k)=\displaystyle\lim_{x\rightarrow \infty} inf \frac{I_{k}(x)}{x},\\ \\
	\overline I_{0}(k)=\displaystyle\lim_{x\rightarrow 0} inf \frac{\overline I_{k}(x)}{x},\hspace{4cm} \overline I_{\infty}(k)=\displaystyle\lim_{x\rightarrow \infty} inf \frac{\overline I_{k}(x)}{x},\\ \\
	I^{q}(k)=\displaystyle\lim_{x\rightarrow q} sup \frac{I_{k}(x)}{x},\ \ \ \hspace{4cm} 	\overline I^{q}(k)=\displaystyle\lim_{x\rightarrow q} sup \frac{\overline I_{k}(x)}{x}.
\end{eqnarray*}

\begin{theorem}
	Assume that the conditions (H1)-(H7) are satisfied. Then the impulsive IBVP (1.1) has at least two positive solutions satisfying $0<\|x_{1}\| <q<\|x_{2}\|$ if for $[a,b]\subset (0,\infty)$, the following conditions hold:
	\begin{itemize}
		\item[$(A1)$] $w\bigg(f_{0} \displaystyle\int_{a}^{b}G(s,s)p(s)ds
		+\frac{\min\{\varphi(\infty),\theta(0)\}}{D}(g_{1_{0}}+g_{2_{0}})\displaystyle\int_{a}^{b}\psi(s)ds \\ \ \ \ \ \ +\displaystyle\sum_{k=1}^{n}G(t_{k},t_{k})\overline I_{0}(k) 
		+\displaystyle\sum_{k=1}^{n}p(t_{k})G_{s}(t,s)|_{\substack{t=t_{k}\\s=t_{k}}}I_{0}(k) \bigg)>1$,\\ \\
		$w\bigg(f_{\infty} \displaystyle\int_{a}^{b}G(s,s)p(s)ds
		+\frac{\min\{\varphi(\infty),\theta(0)\}}{D}(g_{1_{\infty}}+g_{2_{\infty}})\displaystyle\int_{a}^{b}\psi(s)ds \\ +\displaystyle\sum_{k=1}^{n}G(t_{k},t_{k})\overline I_{\infty}(k)
		+\displaystyle\sum_{k=1}^{n}p(t_{k})G_{s}(t,s)|_{\substack{t=t_{k}\\s=t_{k}}}I_{\infty}(k)\bigg)>1$,
		\item[$(A2)$] There exists a $q>0$ such that\\
		$\bigg[1+c\displaystyle\sup_{t\in J}\displaystyle\frac{1}{p(t)}\bigg] \bigg\{h^{q}\displaystyle\int_{0}^{\infty}G(s,s)p(s)k(s)ds
		+\frac{\max\{\varphi(0),\theta(\infty)\}}{D}(g_{1}^{q}+g_{2}^{q}) \\ \times\displaystyle\int_{0}^{\infty}\psi(s)ds
		+\displaystyle\sum_{k=1}^{n}G(t_{k},t_{k})\overline I^{q}(k)
		+\displaystyle\sum_{k=1}^{n}p(t_{k})G_{s}(t,s)|_{\substack{t=t_{k}\\s=t_{k}}}I^{q}(k)\bigg\}<1,\ \ \ \ \ \ \ \\$ for all $\ 0<|x|+|x^{'}| \leq q,\ $ a.e. $t\in [0,\infty).$
	\end{itemize}	
\end{theorem}

\noindent{\em Proof}. By the definition of $f_{0}, I_{0}, \overline I_{0}, g_{1_{0}}$ and $g_{2_{0}}$ for any $\varepsilon>0$, there exist $r\in(0,q)$ such that,
\begin{eqnarray}
\nonumber f(t,x,y) \geq  (1-\varepsilon)f_{0}(|x|+|y|),\ \ \ \ \ (|x|+|y|\leq r,\ \ \ t\in [a,b])
\end{eqnarray}
\begin{eqnarray}\label{c1}
\nonumber g_{1}(x)\geq (1-\varepsilon)g_{1_{0}}x,
\end{eqnarray}
\begin{eqnarray}
g_{2}(x)\geq (1-\varepsilon)g_{2_{0}}x,
\end{eqnarray}
\begin{eqnarray}
\nonumber I_{k}(x)\geq (1-\varepsilon)I_{0}(k)x,
\end{eqnarray}
\begin{eqnarray}
\nonumber \overline I_{k}(x)\geq (1-\varepsilon)\overline I_{0}(k)x,
\end{eqnarray}
\begin{eqnarray}
\nonumber(1-\varepsilon)w\bigg(f_{0} \displaystyle\int_{a}^{b}G(s,s)p(s)ds
+\frac{\min\{\varphi(\infty),\theta(0)\}}{D}(g_{1_{0}}+g_{2_{0}})\displaystyle\int_{a}^{b}\psi(s)ds \\\nonumber 
+\displaystyle\sum_{k=1}^{n}G(t_{k},t_{k})\overline I_{0}(k) 
+\displaystyle\sum_{k=1}^{n}p(t_{k})G_{s}(t,s)|_{\substack{t=t_{k}\\s=t_{k}}}I_{0}(k)\bigg) \geq1.
\end{eqnarray}
Define the open sets
$$\Omega_{r} = \big\{x\in BPC^{1}(J): \|x\|<r \big\}.$$

Let $\Phi=1$ then, $\Phi\in K$. Now we prove that
\begin{eqnarray}\label{c3}
x\neq Tx+\lambda\Phi,\ \ \forall x\in K\cap\partial\Omega_{r}, \ \ \lambda>0 .
\end{eqnarray}

We assume that $x_{0}=Tx_{0}+\lambda_{0}\Phi$ where $x_{0}\in K\cap\partial\Omega_{r}$ and $\lambda_{0}>0$. Let $\ \mu=\displaystyle\min_{t\in[a,b]} x_{0}(t),\ $ then for any $t\in[a,b]$, we have
\begin{eqnarray*}
	& x_{0}(t) & =Tx_{0}(t)+\lambda_{0}\Phi \\&
	& =\displaystyle\int_{0}^{\infty}G(t,s)p(s)f(s,x_{0}(s),x^{'}_{0}(s))ds
	+\frac{\varphi(t)}{D}\displaystyle\int_{0}^{\infty}g_{1}(x_{0}(s))\psi(s)ds
	\\&
	&\ \ \ 
	+\frac{\theta(t)}{D} \displaystyle\int_{0}^{\infty}g_{2}(x_{0}(s))\psi(s)ds +\displaystyle\sum_{k=1}^{n}G(t,t_{k})\overline I_{k}(x_{0}(t_{k}))\\&
	&\ \ \ 
	+\displaystyle\sum_{k=1}^{n}p(t_{k})G_{s}(t,s)|_{s=t_{k}}I_{k}(x_{0}(t_{k}))+\lambda_{0} \\& &
	> w\mu(1-\varepsilon)\bigg\{f_{0}\displaystyle\int_{a}^{b}G(s,s)p(s)ds
	+\frac{\min\{\varphi(\infty),\theta(0)\}}{D}(g_{1_{0}}+g_{2_{0}}) \\&
	&\ \ \ \times\displaystyle\int_{a}^{b}\psi(s)ds +\displaystyle\sum_{k=1}^{n}G(t_{k},t_{k})\overline I_{0}(k)
	+\displaystyle\sum_{k=1}^{n}p(t_{k})G_{s}(t,s)|_{\substack{t=t_{k}\\s=t_{k}}}I_{0}(k)\bigg\} +\lambda_{0} \\& 
	& \geq \mu +\lambda_{0}.
\end{eqnarray*}
This implies $\mu>\mu +\lambda_{0}$ a contradiction. Therefore $\eqref{c3}$ holds.

That the by definition of $f_{\infty}, I_{\infty}, \overline I_{\infty}, g_{1_{\infty}}$ and $g_{2_{\infty}}$ for any $\varepsilon>0$, there exist $R>q$ such that,
\begin{eqnarray}
\nonumber f(t,x,y) \geq(1-\varepsilon)f_{\infty}(|x|+|y|), \ \ \ \ (|x|+|y|\geq R,\ \ \ t\in[a,b]),
\end{eqnarray}
\begin{eqnarray}
\nonumber g_{1}(x)\geq (1-\varepsilon)g_{1_{\infty}}x,\ \ \ \ \ (\forall|x|\geq R),
\end{eqnarray}
\begin{eqnarray}\label{c4}
g_{2}(x)\geq (1-\varepsilon)g_{2_{\infty}}x, 
\end{eqnarray}
\begin{eqnarray}
\nonumber I_{k}(x)\geq (1-\varepsilon)I_{\infty}(k)x,
\end{eqnarray}
\begin{eqnarray}
\nonumber \overline I_{k}(x)\geq (1-\varepsilon)\overline I_{\infty}(k)x,
\end{eqnarray}
\begin{eqnarray}
\nonumber(1-\varepsilon)w\bigg(f_{\infty} \displaystyle\int_{0}^{\infty}G(s,s)p(s)ds
+\frac{\min\{\varphi(\infty),\theta(0)\}}{D}(g_{1_{\infty}}+g_{2_{\infty}})\displaystyle\int_{0}^{\infty}\psi(s)ds \\\nonumber +\displaystyle\sum_{k=1}^{n}G(t_{k},t_{k})\overline I_{\infty}(k) 
+\displaystyle\sum_{k=1}^{n}p(t_{k})G_{s}(t,s)|_{\substack{t=t_{k}\\s=t_{k}}}I_{\infty}(k)\bigg) \geq1.
\end{eqnarray}
Define the open sets
\begin{eqnarray}\label{c5}
\Omega_{R} = \big\{x\in BPC^{1}(J): \|x\|<R \big\}.
\end{eqnarray}
As the proof of $\eqref{c3}$, we can get that
\begin{eqnarray}\label{c6}
x\neq Tx+\lambda\Phi,\ \ \forall x\in K\cap\partial\Omega_{R},\ \ \ \lambda>0.
\end{eqnarray}
On the other hand, for any $\varepsilon>0$, choose $q$ in $(A2)$ such that
\begin{eqnarray}\label{c7}
\nonumber (1+\varepsilon)\bigg[1+c\sup_{t\in J}\displaystyle\frac{1}{p(t)}\bigg] \bigg\{ h^{q}\displaystyle\int_{0}^{\infty}G(s,s)p(s)k(s)ds
+\frac{\max\{\varphi(0),\theta(\infty)\}}{D}(g_{1}^{q}+g_{2}^{q})\\ \\\nonumber \times \displaystyle\int_{0}^{\infty}\psi(s)ds
+\displaystyle\sum_{k=1}^{n}G(t_{k},t_{k})\overline I^{q}(k)
+\displaystyle\sum_{k=1}^{n}p(t_{k})G_{s}(t,s)|_{\substack{t=t_{k}\\s=t_{k}}}I^{q}(k) \bigg\}\leq 1.
\end{eqnarray}

By the definition of  $\ \ h^{q},\ I^{q},\ \overline I^{q},\ g_{1}^{q}$ and $\ g_{2}^{q},\ $ for the above $\varepsilon>0$, there exists $\delta>0$, when $\ \ |x|, |x|+|y|\in (q-\delta,q+\delta)$; thus, we get
\begin{eqnarray}\label{c8}
\nonumber h(x,y) \leq(1+\varepsilon)h^{q}(|x|+|y|),\\\nonumber\\ 
\nonumber g_{i}(x) \leq(1+\varepsilon)g_{i}^{q}x,\ \ \ 1\leq i\leq2 \\ \\
\nonumber I_{k}(x) \leq(1+\varepsilon)I^{q}(k)x, \\\nonumber\\
\nonumber \overline I_{k}(x)\leq (1+\varepsilon)\overline{I}^{q}(k)x.
\end{eqnarray}
Define the open sets
\begin{eqnarray}\label{c9}
\Omega_{q} = \big\{x\in BPC^{1}(J): \|x\|<q \big\}.
\end{eqnarray}
Then for any $x\in K\cap\partial\Omega_{q}$ and $t\in J$ we obtain that,
\begin{eqnarray}
& \nonumber |(Tx)(t)|+&|(Tx)'(t)|\leq\displaystyle\int_{0}^{\infty}G(s,s)p(s)f(s,x(s),x^{'}(s))ds\\&\nonumber&\ \ \ \ +\displaystyle\frac{\max\{\varphi(0),\theta(\infty)\}}{D}\displaystyle\int_{0}^{\infty}[g_{1}(x(s))+g_{2}(x(s))]\psi(s)ds\\&\nonumber&\ \ \ \ +\displaystyle\sum_{k=1}^{n}G(t_{k},t_{k})\overline I_{k}(x(t_{k}))\\&\nonumber&\ \ \ \ 
+\displaystyle\sum_{k=1}^{n}p(t_{k})G_{s}(t,s)|_{\substack{t=t_{k}\\s=t_{k}}}I_{k}(x(t_{k}))\\&\nonumber&\ \ \ \ 
+c\displaystyle\sup_{t\in J} \displaystyle\frac{1}{p(t)}\displaystyle\int_{0}^{\infty}G(s,s)p(s)f(s,x(s),y(s))ds\\&\nonumber&\ \ \ \ +\displaystyle\frac{\max\{a_{1},a_{2}\}}{D}\displaystyle\sup_{t\in J}\displaystyle\frac{1}{p(t)}\displaystyle\int_{0}^{\infty}[g_{1}(x(s))+g_{2}(x(s))]\psi(s)ds\\&\nonumber&\ \ \ \ 
+c\displaystyle\sup_{t\in J} \displaystyle\frac{1}{p(t)}[ \displaystyle\sum_{k=1}^{n}G(t_{k},t_{k})\overline I_{k}(x(t_{k})) \\&\nonumber&\ \ \ \ +\displaystyle\sum_{k=1}^{n}p(t_{k})G_{s}(t,s)|_{\substack{t=t_{k}\\s=t_{k}}}I_{k}(x(t_{k}))]
\end{eqnarray}

\begin{eqnarray}
&\nonumber& \leq\bigg[1+c\displaystyle\sup_{t\in J} \displaystyle\frac{1}{p(t)}\bigg]\displaystyle\int_{0}^{\infty}G(s,s)p(s)k(s)h(x(s),x^{'}(s))ds\\&\nonumber&\ \ \ \ +\bigg[\displaystyle\frac{\max\{\varphi(0),\theta(\infty)\}}{D}+\displaystyle\frac{\max\{a_{1},a_{2}\}}{D}\displaystyle\sup_{t\in J} \displaystyle\frac{1}{p(t)}\bigg]\\&\nonumber&\ \ \ \ \times \displaystyle\int_{0}^{\infty}[g_{1}(x(s))+g_{2}(x(s))]\psi(s)ds\\&\nonumber&\ \ \ \ +\bigg[1+c\displaystyle\sup_{t\in J} \displaystyle\frac{1}{p(t)}\bigg]\bigg[ \displaystyle\sum_{k=1}^{n}G(t_{k},t_{k})\overline I_{k}(x(t_{k}))\\&\nonumber&\ \ \ \
+\displaystyle\sum_{k=1}^{n}p(t_{k})G_{s}(t,s)|_{\substack{t=t_{k}\\s=t_{k}}}I_{k}(x(t_{k}))\bigg]\\&\nonumber& \leq (1+\varepsilon)\bigg[1+c\displaystyle\sup_{t\in J} \displaystyle\frac{1}{p(t)}\bigg]\bigg[h^{q}\displaystyle\int_{0}^{\infty}G(s,s)p(s)k(s)ds+\displaystyle\frac{\max\{\varphi(0),\theta(\infty)\}}{D}\\&\nonumber&\ \ \ \ \times(g_{1}^{q}+g_{2}^{q})\displaystyle\int_{0}^{\infty}\psi(s)ds +\displaystyle\sum_{k=1}^{n}G(t_{k},t_{k})\overline I^{q}(k) \\&\nonumber&\ \ \ \
+\displaystyle\sum_{k=1}^{n}p(t_{k})G_{s}(t,s)|_{\substack{t=t_{k}\\s=t_{k}}}I^{q}(k)\bigg] \|x\| \\& \nonumber & \leq \|x\|.
\end{eqnarray}
Therefore $\|Tx\|\leq\|x\|$.

So, we can get the existence of two positive solutions $x_{1}$ and $x_{2}$ satisfying $0<\|x_{1}\|<q<\|x_{2}\|$ by using lemma $\eqref{b4}$.

Using a similar proof of Theorem $(3.1)$, we can get the following theorem.
\begin{theorem}
	Assume that the conditions (H1)-(H7) are satisfied. Then the impulsive IBVP $\eqref{a1}$ has at least two positive solutions satisfying $0<\|x_{1}\| <q<\|x_{2}\|$ if for $[a,b]\subset (0,\infty)$, the following conditions hold;
	
	\begin{itemize}
		\item[$(A3)$]  $\bigg[1+c\displaystyle\sup_{t\in J} \displaystyle\frac{1}{p(t)}\bigg]\bigg\{ h^{0}\displaystyle\int_{0}^{\infty}G(s,s)p(s)k(s)ds
		+\frac{\max\{\varphi(0),\theta(\infty)\}}{D}(g_{1}^{0}+g_{2}^{0}) \\ \ \ \ \ \times\displaystyle\int_{0}^{\infty}\psi(s)ds +\displaystyle\sum_{k=1}^{n}G(t_{k},t_{k})\overline{I}^{0}(k) 
		+\displaystyle\sum_{k=1}^{n}p(t_{k})G_{s}(t,s)|_{\substack{t=t_{k}\\s=t_{k}}}I^{0}(k) \bigg\}<1$,\\ \\
		$\bigg[1+c\displaystyle\sup_{t\in J} \displaystyle\frac{1}{p(t)}\bigg]\bigg\{h^{\infty} \displaystyle\int_{0}^{\infty}G(s,s)p(s)k(s)ds
		+\frac{\max\{\varphi(0),\theta(\infty)\}}{D}(g_{1}^{\infty}+g_{2}^{\infty}) \\ \ \ \ \  \times\displaystyle\int_{0}^{\infty}\psi(s)ds +\displaystyle\sum_{k=1}^{n}G(t_{k},t_{k})\overline{I}^{\infty}(k) 
		+\displaystyle\sum_{k=1}^{n}p(t_{k})G_{s}(t,s)|_{\substack{t=t_{k}\\s=t_{k}}}I^{\infty}(k)\bigg\}<1$,
		\item[$(A4)$] There exists a $q>0$ such that\\
		$w\bigg\{f_{q}\displaystyle\int_{a}^{b}G(s,s)p(s)ds
		+\frac{\min\{\varphi(\infty),\theta(0)\}}{D}(g_{1_{q}}+g_{2_{q}}) \displaystyle\int_{a}^{b}\psi(s)ds\\ \ \ \ \ +\displaystyle\sum_{k=1}^{n}G(t_{k},t_{k})\overline{I}_{q}(k)
		+\displaystyle\sum_{k=1}^{n}p(t_{k})G_{s}(t,s)|_{\substack{t=t_{k}\\s=t_{k}}}I_{q}(k)\bigg\}>1,\ \ \ \ \ \ \\ $ for all $\ 0<|x|+|x^{'}| \leq q,\ $ a.e. $t\in [0,\infty). $
	\end{itemize}	
\end{theorem}

\section{Example}
To illustrate how our main results can be used in practise we present the folloing example.

Consider the following boundary value problem:
\begin{eqnarray}\label{d1}
\begin{cases}
e^{-t}(e^{t}x'(t))'+f(t,x(t),x^{'}(t))=0,\ \ t\in J_{+},\ t\neq \frac{1}{2},\\ \\

\Delta x |_{t=\frac{1}{2}}=\displaystyle\frac{1}{100}x(\frac{1}{2})\\ \\

\Delta x'|_{t=\frac{1}{2}}=\displaystyle\frac{1}{2-e^{-\frac{1}{2}}}x^{-1}(\frac{1}{2}) \\ \\

x(0)=\displaystyle\frac{1}{40\pi}\displaystyle\int_{0}^{\infty} \displaystyle\frac{x^{2}(s)}{1+s^{2}}ds,\\ \\

\displaystyle\lim_{t\rightarrow \infty}e^{t}x'(t)=\frac{1}{8000\pi}\displaystyle\int_{0}^{\infty}\displaystyle\frac{x^{2}(s)}{1+s^{2}} ds,\\ 
\end{cases} 
\end{eqnarray}
where $f(t,x(t),x^{'}(t))=\displaystyle\frac{e^{-t}}{e^{t}-2},\ a_{1}=1, a_{2}=0,\  b_{1}=1,\ b_{2}=1,\ p(t)=e^{t},\ \psi(t)=\displaystyle\frac{1}{1+t^{2}},\\ \\ I_{k}(x(t))=\displaystyle\frac{x(t)}{100},\ \overline I_{k}(x(t))=\displaystyle\frac{x^{-1}(t)}{2-e^{-\frac{1}{2}}},\ g_{1}(x(s))=\frac{x^{2}(s)}{40\pi},\ g_{2}(x(s))=\frac{x^{4}(s)}{8000\pi}$.\\

Set $k(t)=\displaystyle\frac{e^{-t}}{e^{t}-2}$,
$h(x(t),y(t))=1 $ and $q=10$.
It follows from a direct calculation that\\ \\ $g_{1}^{q}=\displaystyle\frac{1}{4\pi},\ g_{2}^{q}=\displaystyle\frac{1}{8\pi},\ g_{1_{\infty}}=\infty,\ g_{2_{\infty}}=\infty,\ g_{1_{0}}=0,\ g_{2_{0}}=0,\ I^{q}(k)=\displaystyle\frac{1}{100},\ \overline I^{q}(k)=\displaystyle\frac{1}{100(1-e^{-\frac{1}{2}})},\\ I_{0}(k)=\displaystyle\frac{1}{100},\ \overline I_{0}(k)=\infty,\ I_{\infty}(k)=\displaystyle\frac{1}{100},\ \overline I_{\infty}(k)=0.$\\ \\
Furthermore, $f_{0}=\infty,\ f_{\infty}=1,$ $\displaystyle\int_{0}^{\infty}G(s,s)p(s)k(s)ds=1$ and  $\displaystyle\int_{0}^{\infty}\psi(s)ds=\displaystyle\frac{\pi}{2}.$ Thus $(A1)$ and $(A2)$ are satisfied. Therefore, by Theorem $\eqref{c1}$, the impulsive IBVP $\eqref{d1}$ has at least two positive solutions $x_{1},x_{2}$ satisfying $0<\|x_{1}\|<10<\|x_{2}\|. $

\end{document}